\documentclass[a4paper,11pt]{article}
\usepackage{latexsym}
\usepackage{amsfonts,amssymb}
\usepackage{amscd}
\usepackage[dvips]{graphicx}

\newtheorem{defn}{Definition}

\newtheorem{observen}{Observation}
\newtheorem{lemman}{Lemma}

\begin{document}
\begin{center}
\noindent {\bf\Large On GCD-morphic sequences}

\vspace{0.5cm}

\noindent {M. Dziemia\'nczuk(*), Wies\l aw Bajguz(**)}

\vspace{0.5cm}

\noindent {\small Institute of Computer Science, Bia\l ystok University \\
PL-15-887 Bia\l ystok, ul.Sosnowa 64, POLAND \\
e-mail: (*) Maciek.ciupa@gmail.com, (**) Bajguz@wp.pl} \\
\end{center}

\begin{abstract}
	This note is a response to one of the problems posed by Kwa\'sniewski in \cite{1,2}, see also \cite{3} i.e. GCD-morphic Problem III. 
We show that any GCD-morphic sequence $F$ is at the point product of primary GCD-morphic sequences and any GCD-morphic sequence is encoded by natural number valued sequence satisfying condition (C1). The problem of general importance - for example in number theory was formulated in \cite{1,2}
while investigating a new class of DAG's and their correspondent p.o. sets encoded uniquely by sequences with combinatorially interpretable properties.

\end{abstract}

\vspace{0.4cm}

\noindent Key Words: GCD-morphic sequence

\vspace{0.1cm}

\noindent AMS Classification Numbers: 06A07, 05A10, 11A41, 05C20
\vspace{0.1cm}

\noindent  Presented at Gian-Carlo Polish Seminar:\\
\noindent \emph{http://ii.uwb.edu.pl/akk/sem/sem\_rota.htm}

\vspace{0.4cm}

\section{Preliminaries}

\begin{defn} [GCD-morphic sequence \cite{1}]
The integers' valued  sequence $F=\{F_n\}_{n\geq 0}$  is called the
GCD-morphic  sequence  if  $$GCD(F_n,F_m) = F_{GCD(n,m)}$$
where $GCD$ stays for Greatest Common Divisor operator.
\end{defn}

\noindent \textbf{GCD-morphism Problem. Problem III. \cite{1}} \textit{ Find effective characterizations
and/or  an algorithm to produce the GCD-morphic sequences i.e. find all examples.}

\begin{defn}
Any natural numbers valued sequence $G$ of the form
$$
	G_{c,N} \equiv \{g_n\}_{n\geq 0},\ where \ g_n = \left\{ 
	\begin{array}{ll}
		c & N | n \\ 
		1 & n=0 \vee N\nmid{n}
		\end{array}
	\right.
$$
is called the primary GCD-morphic sequence. 
\end{defn}

\begin{observen}
For any natural $c>0$, $G_{c,k}$ is GCD-morphic.
\end{observen}

\noindent \textbf{Notation}

\noindent The product of two sequences $F_n=(F_1^n,F_2^n,\dots)$, $F_k=(F_1^k,F_2^k,\dots)$ is at the point product of functions i.e. 
$F_n\cdot F_k=(F_1^n\cdot F_1^k,F_2^n\cdot F_2^k,\dots)$.

\bigskip

\section {GCD-morphic sequences}

\begin{lemman} 
Let $F_n=( 1,\dots,1,F^n_n,F^n_{n+1},\dots)$ be a GCD-morphic sequence. Then we have the equality $$F_n=G_{F_n^n,n}\cdot F_{n+1} ,$$ 
where \\ 
$F_{n+1}=\lbrace F_k^{n+1}\rbrace _{k=1}^\infty $, \ 
$F_k^{n+1}=
  \cases{
     1,       & if $k\le n$;   \cr 
     {F_k^n}\over{F_n^n}, & if $n\vert k$; \cr 
     F_k^n, & otherwise. \cr
     }
$ \\
and $F_{n+1}$ is GCD-morphic.

\end{lemman}

\noindent \textbf{Proof.}

Since $F_n$ is GCD-morphic, it follows that  $F_n^n\vert F_k^n$ for $n\vert k$. Therefore $F_k^{n+1}$ is well defined and the equality is obvious.

Now we must proof that $F_{n+1}$ is GCD-morphic. \\
\noindent Let $k,l$ be any natural numbers. We must show, that $CGD$-morphic condition: $GCD(F_k^{n+1},F_l^{n+1})=F_{GCD(k,l)}^{n+1}$ is true. In order this to be proved  let us investigate the following cases:

a) $n\not\vert k$ and $n\not\vert l$,

b) $n\vert k$ and $n\vert l$,

c) otherwise.

 In the  a) case $F_k^{n+1}=F_k^n$, $F_l^{n+1}=F_l^n$ and $F_{GCD(k,l)}^{n+1}=F_{GCD(k,l)}^n$ holds. Since $F_n$ is GCD-morphic, therefore the $CGD$-morphic condition is satified.

 In the case b) we obtain $n\vert GCD(k,l)$ and accordingly \\ 
$F_k^{n+1}={F_k^n \over F_n^n}$, $F_l^{n+1}={F_l^n \over F_n^n}$ and $F_{GCD(k,l)}^{n+1}={F_{GCD(k,l)}^n \over F_n^n}$. \\
It provides to $GCD(F_k^{n+1},F_l^{n+1})={GCD(F_k^n,F_l^n) \over F_n^n}=
{F_{GCD(k,l)}^n \over F_n^n}=F_{GCD(k,l)}^{n+1}$

 To prove that  $CGD$-morphic condition is satified in the case c) let us assume that 
$n\vert k$ and $n\not\vert l$. Then \\
$F_k^n=F_k^{n+1} \cdot F_n^n$,  $F_l^{n+1}=F_l^n$ and  $F_{GCD(k,l)}^{n+1}=F_{GCD(k,l)}^n$
and consequently\\
$GCD(F_k^{n+1} \cdot F_n^n , F_l^{n+1})=GCD(F_k^n , F_l^n) = 
F_{GCD(k,l)}^n=F_{GCD(k,l)}^{n+1}$.

\noindent A contrario proof: assume that 
$GCD(F_k^{n+1} \cdot F_n^n , F_l^{n+1}) \ne GCD(F_k^{n+1} , F_l^{n+1})$. It is equivalent to $GCD(F_k^{n+1} \cdot F_n^n , F_l^n) \ne GCD(F_k^{n+1} , F_l^n)$.\\
Then $GCD(F_n^n,F_l^n)\ne 1$, but otherwise $GCD(n,l)<n$ and consequently $GCD(F_n^n,F_l^n)=F_{GCD(n,l)}^n=1$. Hence - contradiction.

\noindent Therefore $GCD(F_k^{n+1} , F_l^{n+1})=GCD(F_k^{n+1} \cdot F_n^n , F_l^{n+1})=F_{GCD(k,l)}^{n+1}$.$\blacksquare$ 

\bigskip
Since every periodic sequence $((1,\dots,1,F^n_n))$ is GCD-morphic then according to the above lemma we sum up.
Lemma 1 provides a prescription for how  to produce infinite number of primary GCD-morphic sequences. Recall:  every such product sequence is well defined i.e. each element of this sequence is a product of \textbf{finite} number elements different from 1.

\medskip

\noindent \textbf{Conclusion}

\noindent Let $F=(F^1_1,F^1_2,\dots)$ be a GCD-morphic sequence. Then $F$ is a product $G_{c_1,1} \cdot G_{c_2,2} \cdot\dots$ of primary GCD-morphic sequences, where 
$$c_n={F_n^1\over \prod_{k\vert n, k<n}c_k} \ .$$

\medskip

\begin{lemman}
	Let a sequence $F=(F_1,F_2,\dots )$ be given as the product of infinite number primary GCD-morphic sequences $G_{c_i,i}$, $i\in N$. The sequence $F$ is GCD-morphic if, and only if for any $k,n\in N$\\

\noindent $(C1) \ \ \ \ \ \ \ \ \ \ \ \ \ k<n \wedge GCD(n,k)\neq k \Rightarrow GCD(c_n, c_k) = 1$

\end{lemman}

\noindent \textbf{Proof.}

\noindent Let any GCD-morphic sequence F be given.\\ 
To a contrario proof assume that $(C1)$ \textbf{is not} true. Consider such natural numbers $k,n$ that\\ 
(*) \ \ $k<n$,  $GCD(n,k)=\alpha \neq k$ and $GCD(c_n, c_k) > 1$.\\ 
Then $GCD(F_n, F_k) = F_\alpha$. \\
Since $GCD(n,k)=\alpha \neq k$, then $\alpha | k \wedge \alpha < k$, and therefore\\ 
(1) \ \ $F_k = \prod_{j|k}{c_j} = \prod_{j|\alpha}\cdot c_k \cdot F_k'$. \\
Moreover from $GCD(n,k)=\alpha \neq k$ we have $\alpha | n \wedge \alpha < n$ and it leads to\\
(2) \ \ $F_n = \prod_{j|n}{c_j} = \prod_{j|\alpha}\cdot c_n \cdot F_n'$. \\ 
From (1) and (2) we have that \\
$GCD(F_n, F_k) = \prod_{j|\alpha}{c_j} \cdot GCD(c_k\cdot F_k', c_n\cdot F_n') \geq F_\alpha \cdot GCD(c_k, c_n)$\\
By (*) $GCD(c_n, c_k) > 1$, and therefoe  $GCD(F_n, F_k) > F_\alpha$. Hence contradiction.

\vspace{0.4cm}
\noindent Next, we need to show that if $F$ is a product of infinite number primary GCD-morphic sequences $G_{c_i,i}$ which fulfils (C1) then $F$ is GCD-morphic. \\
Since $F$ is a product of $G_{c_i,i}$, then $F_n=\prod_{j|n}{c_j}$ for any $n$.\\ 
Therefore \ \ 
$GCD(F_n, F_k) = GCD(\prod_{j|n}{c_j}, \prod_{j|k}{c_j}) =\\
= GCD(\prod_{j|n \wedge j|k}{c_j}\cdot\prod_{j|n\wedge j \nmid k}{c_j}, \prod_{j|n \wedge j|k}{c_j}\cdot\prod_{j|k\wedge j \nmid n}{c_j}) = \prod_{j|n \wedge j|k}{c_j}\cdot P = \\
= \prod_{j|GCD(n,k)}{c_j}\cdot P = F_{GCD(n,k)} \cdot P$ \\ 
where $P = GCD(\prod_{j|n \wedge j\nmid k}{c_j}, \prod_{j|k \wedge j\nmid n}{c_j})$. \\
If $P > 1$ then there exists r,s such that\\
(*) \ \ \ $r|n, r\nmid k, s|k, s\nmid n$ and $GCD(c_r, c_s) > 1$. \\
By (*) holds $r\nmid s$ and $s\nmid r$. Then $GCD(c_r, c_s) = 1$ by (C1) - contrary to (*). Hence $P=1$ and $GCD(F_n, F_k)= F_{GCD(n,k)}$. It ends the proof $\blacksquare$

\vspace{0.4cm}

The lemmas 1 and 2 lead to the following conclusion.

\vspace{0.4cm}

\noindent \textbf{Conclusion}

\noindent Any sequence $C\equiv\{c_n\}_{n\geq 1}$ such that $\{ G_{c_n,n}\}_{n\geq 1}$ satisfies condition $(C1)$ encodes GCD-morphic sequence and otherwise - any GCD-morphic sequence is encoded by such sequence. The correspondence is biunivoque.


\vspace{0.4cm}

\noindent \textbf{Examples:}

\noindent Natural numbers' sequence $C = (1,2,3,2,5,1,7,2,3,1,11,1,13,1,1,2,17,...)$

\noindent Fibonacci numbers' sequence $C = (1,1,2,3,5,4,13,7,17,11,89,6,...)$

\noindent The sequence of products primary divisors of natural numbers and it is a GCD-morphic sequence encoded by
 $c_n=  \cases{
     n       & if $n$ primary;   \cr 
     1 & otherwise. \cr
     }$. 

\bigskip

 \noindent \textbf{\Large Acknowledgements}

\medskip

We  would like to thank Professor A. Krzysztof Kwa\'sniewski for guidance and final improvements of this paper. Attandence
of Ewa Krot-Sieniawska is highly appreciated too.

\vspace{0.2cm}
 

\end{document}